\documentclass[12pt]{article}
\usepackage{url}
\usepackage{amsmath, amssymb}

\usepackage[margin=30mm]{geometry}

\newcommand{\of}{\Omega_{\mathcal F}}
\newcommand{\hof}{{\hat{\Omega}}_{\mathcal F}}

\renewcommand{\hat}{\widehat}
\newcommand{\Gd}{G_{\diamond}}
\newcommand{\Vd}{V_{\diamond}}
\newcommand{\Ed}{E_{\diamond}}
\renewcommand{\c}{\widetilde{c}}
\newcommand{\Z}{\mathbb{Z}}
\newcommand{\R}{\mathbb{R}}
\newcommand{\Ztaxi}{\vec{\mathbb Z}^2}

\renewcommand{\Gd}{G_{\Diamond}}
\renewcommand{\Vd}{V_{\Diamond}}
\renewcommand{\Ed}{E_{\Diamond}}
\newcommand{\Ge}{G_0}
\newcommand{\Go}{G_1}

\newcommand{\Gb}{G_b}
\newcommand{\Eb}{E_b}
\newcommand{\Vb}{V_b}
\newcommand{\hatG}{\hat{G}}

\newcommand{\qed}[0]{{\hspace*{\fill}\mbox{$\Box$}}}

\newcommand{\besttorus}{5.3646}
\newcommand{\bestbox}{7.1031}
\newcommand{\bestmuupper}{1.5884}
\newcommand{\bestmulower}{1.5196}
\newcommand{\bestlambdalower}{4.3332}

\newtheorem{theorem}{Theorem}[section]
\newtheorem{definition}{Definition}[section]
\newtheorem{lemma}[theorem]{Lemma}
\newtheorem{corollary}[theorem]{Corollary}

\makeatletter
\def\bigpar{\bigbreak\@afterindentfalse\@afterheading\ignorespaces}
\def\medpar{\medbreak\@afterindentfalse\@afterheading\ignorespaces}
\def\smallpar{\smallbreak\@afterindentfalse\@afterheading\ignorespaces}

\makeatletter
\def\eqalign#1{\,\vcenter{\openup\jot\m@th
  \ialign{\strut\hfil$\displaystyle{##}$&$\displaystyle{{}##}$\hfil
     \crcr#1\crcr}}\,}
\makeatother

\begin{document}

\title{Phase Coexistence and Slow Mixing for the Hard-Core Model on $\Z^2$}
\author{Antonio Blanca\thanks{Computer Science Division, U.C. Berkeley, Berkeley, CA 94720, tonyblanca@gmail.com}, \ \ David Galvin\thanks{Department of Mathematics, University of Notre Dame, Notre Dame, IN 46656, dgalvin1@nd.edu}, \ \ Dana Randall\thanks{School of Computer Science, Georgia Institute of Technology, Atlanta, GA 30332-0765, randall@cc.gatech.edu} \ \ and \  Prasad Tetali\thanks{School of Mathematics, Georgia Institute of Technology, Atlanta, GA 30332-0280, tetali@math.gatech.edu}}
\date{\today}

\maketitle

\begin{abstract}
The hard-core model has attracted much attention across several disciplines, representing lattice gases in statistical physics and independent sets in discrete mathematics and computer science.
On finite graphs, we are given a parameter $\lambda$, and an independent set $I$ arises with probability proportional to $\lambda^{|I|}$.  On infinite graphs a Gibbs distribution is defined 
as a suitable limit with the correct conditional probabilities.
In the infinite setting we are interested in determining when this limit is unique and when there is phase coexistence, i.e., existence of multiple Gibbs states.
On finite graphs we are interested in determining the mixing time of local Markov chains.

On $\Z^2$  it is conjectured that these problems are related and that both undergo a phase transition at some critical point $\lambda_c \approx 3.79$ \cite{approx}. For the question of phase coexistence much of the work to date has focused on the regime of uniqueness, with the best result to date being recent work of Restrepo et al. \cite{rstvy} showing that there is a unique Gibbs state for all $\lambda < 2.3882$. Here we give the first non-trivial result in the other direction, showing that there are multiple Gibbs states for all $\lambda > \besttorus$. Our proof begins along the lines of the standard Peierls argument, but we add two significant innovations. First, building on the idea of fault lines introduced by Randall \cite{randall}, we construct an event that distinguishes two boundary conditions and yet always has long contours associated with it, obviating the need to accurately enumerate short contours. Second, we obtain vastly improved bounds on the number of contours by relating them to a new class of self-avoiding walks on an oriented version of $\Z^2$.

The best  result for rapid mixing of local Markov chains on boxes of $\Z^2$ is also when $\lambda < 2.3882$ \cite{rstvy}.
Here we extend our characterization of fault lines to show that local Markov chains will mix slowly when $\lambda > \besttorus$ on lattice regions with periodic (toroidal) boundary conditions and when $\lambda > \bestbox$ with non-periodic (free) boundary conditions. The arguments here rely on a careful analysis that relates contours to taxi walks and represent a sevenfold improvement to the previously best known values of $\lambda$ \cite{randall}.
\end{abstract}

\section{Introduction}
The hard-core model was introduced in statistical physics as a model for lattice gases, where each molecule occupies non-trivial space in the lattice, requiring occupied sites to be non-adjacent.
When a 
lattice such as $\Z^d$ is viewed as a graph, the allowed configurations of molecules naturally correspond to  independent sets in the graph.

Given a graph $G$, let ${\Omega}$ be the set of independent sets of~$G$.  Given a (fixed) {\it activity} (or {\it fugacity}) $\lambda \in \R^+$, the weight associated with each independent set $I$ is $w(I) = \lambda^{|I|}$.  The associated {\it Gibbs} (or {\it Boltzmann}) {\it distribution} $\mu = \mu_{G, \lambda}$ is defined on ${\Omega}$, assuming $G$ is finite, as $\mu(I) = w(I)/Z$, where the normalizing constant $Z = Z(G, \lambda) = \sum_{J \in {\Omega}} w(J)$ is commonly called the {\it partition function}.
Physicists are interested in the behavior of models on an infinite graph (such as the integer lattice $\Z^d$), where the Gibbs measure is defined as a certain weak limit with appropriate conditional probabilities.  For many models  it is believed that as a parameter of the system is varied --  such as the inverse temperature $\beta$ for the Ising model or the activity $\lambda$ for the hard-core model --  the system undergoes a phase transition at a critical point.

For the classical Ising model,
Onsager, in seminal work \cite{onsager},
established the precise value of the critical temperature $\beta_c(\Z^2)$ to be $\log(1+\sqrt{2})$.  Only recently have the analogous values for the (more general) $q$-state Potts model been established in breakthrough work by Beffara and Duminil-Copin \cite{bd-c}, settling a more than half-a-century old open problem. Establishing such a precise value for the hard-core model with
currently available methods seems nearly impossible.
Even the {\em existence} of such a (unique) critical activity $\lambda_c$, where there is a transition from a unique Gibbs state to the coexistence of multiple Gibbs states remains conjectural for $\Z^d$ ($d\ge 2$; it is folklore that there is no such transition for $d=1$), while it is simply  untrue for
general graphs
(even general trees, in fact, thanks to a result of Brightwell et al.~\cite{bhw}).  Regardless,  a non-rigorous prediction from the statistical physics literature \cite{approx} suggests $\lambda_c \approx 3.796$ for $\Z^2$.

Thus, from a statistical physics or  probability point of view, understanding the precise dependence on $\lambda$ for the existence of unique or multiple hard-core Gibbs states is a natural and challenging problem.
Moreover, breakthrough works of Weitz \cite{weitz}  and Sly \cite{sly} in recent years identified $\lambda_c({\mathbb T}_\Delta)$ -- the critical activity for the hard-core model on an infinite $\Delta$-regular tree --  as a {\em computational} threshold where estimating the hard-core partition function on general $\Delta$-regular graphs undergoes a transition from being in $P$ to being $NP$-hard (specifically, there is no PTAS unless $NP = RP$), further motivating the study of such (theoretical) physical transitions and their computational implications. While it is not surprising that for many fundamental problems computing the partition function {\em exactly} is intractable, it is remarkable that even approximating the partition function of the hard-core model above a certain critical threshold also turns out to be hard.

Starting with Dobrushin \cite{dob}  in 1968, physicists
have been developing techniques
to 
characterize the regimes on either side of $\lambda_c$ for the hard-core model. Most attention has focused on establishing ever larger values of $\lambda$ below which there is always uniqueness of phase. The problem has proved to be a fruitful one for the blending of ideas from physics, discrete probability and theoretical computer science, with improvements to our understanding of the problem having been made successively by Radulescu and Styer \cite{radulescustyer}, van den Berg and Steif \cite{bs} and Weitz \cite{weitz}, among others. The state of the art is recent work of Restrepo et al. \cite{rstvy}, building on the novel arguments introduced by Weitz, which establishes uniqueness for all $\lambda < 2.3882$.

Much less is known about the regime of phase coexistence. Dobrushin \cite{dob} established phase coexistence for all $\lambda > C$, but did not explicitly calculate $C$. Around the time of the writing of \cite{bc7}, Borgs reported that $C=80$ was the {\em theoretical} lower limit of Dobrushin's argument \cite{borgs_pc}, but a recent computation by the second author suggested that the {\em actual} consequence of the argument was more like $C \approx 300$.

\medskip

\noindent {\bf Our first main result}: Using insights from physics and combinatorics, we establish the first non-trivial upper bound on $\lambda_c$ for the hard-core model on $\Z^2$.
\begin{theorem} \label{thm-phasecoexistence}
For all $\lambda > 5.3646$, the hard-core model on $\Z^2$ with activity $\lambda$ admits multiple Gibbs states.
\end{theorem}

\medskip

From a computational standpoint, there are two natural questions to ask concerning the hard-core model on a finite graph. Can the partition function be approximated? (an exact determination is in general too much to hope for), and how easy is it to sample from a given Gibbs distribution? For both questions, a natural and powerful method is offered by Markov chain algorithms -- carefully constructed random walks on the space of independent sets of a graph whose equilibrium distributions are the desired Gibbs distributions. One of the most commonly studied families of Markov chains are the {\em local-update} chains, such as Glauber dynamics, that change the state of a bounded number of vertices at each step (in particular, Glauber dynamics changes the state of a single vertex at each step).

The efficiency of the Markov chain method relies on the underlying chain being rapidly mixing; that is, it must fairly quickly reach a distribution that is quantifiably close to stationary.
For many problems, local chains such as Glauber dynamics seem to mix rapidly below some critical point, while mixing slowly above that point.  Most notably for the Ising model on $\Z^2$,  simple local Markov chains are rapidly mixing (in fact, with optimal rate) for $\beta < \beta_c(\Z^2)$ and slowly mixing
for $\beta > \beta_c(\Z^2)$; various mathematical physics experts (including Aizenmann, Holley, Martinelli, Olivieri, Schonmann, Stroock and Zegarlinski) have contributed to this work. Recently the Ising picture was completed by Lubezky and Sly \cite{LubezkySly}, who showed polynomial mixing at $\beta = \beta_c(\Z^2)$.

Once again, the known bounds are less sharp for the hard-core model. Luby and Vigoda \cite{lv} showed that Glauber dynamics on independent sets is fast when  $\lambda \leq 1$ on the $2$-dimensional lattice and torus.   Weitz \cite{weitz} reduced the analysis on the grid to the tree, thus establishing that in this same setting Glauber dynamics is fast up to the critical point for the $4$-regular tree, in effect for $\lambda < 1.6875$. The state of the art in this direction is due to Restrepo et al. \cite{rstvy}. Using now standard machinery (establishing so-called {\em strong spatial mixing}), they proved that the natural Glauber dynamics on the space of hard-core configurations on boxes in $\Z^2$ is rapidly mixing for all $\lambda < 2.3882$. These results also lead to efficient deterministic algorithms for approximating the hard-core partition function
on (finite regions of) $\Z^2$.

As with phase-coexistence, it is believed that there is a critical value $\lambda_c^{\rm mix}$ across which the Glauber dynamics for sampling from hard-core configurations on a $n$ by $n$ box in $\Z^2$ flips from mixing in time polynomial in $n$, to exponential in $n$, and that it coincides with $\lambda_c$.
Pinning down such a precise transition seems beyond reach at the moment, but we can try to establish a short range of values into which $\lambda_c^{\rm mix}$ (if it exists) must fall, by finding ever smaller bounds for $\lambda$ above which the mixing time is exponential.

Borgs et al. \cite{bc7} showed that Glauber dynamics is slow
on toroidal lattice regions in $\Z^d$ (for $d\ge 2$),
when $\lambda$ is  sufficiently large,
in particular establishing a finite constant above which mixing is slow on $\Z^2$. The first effective bound was provided by Randall  \cite{randall}, who showed slow mixing for $\lambda > 50.59$ on boxes with periodical boundary conditions, and for $\lambda > 56.88$ on boxes with free boundary. (Better bounds were originally reported in \cite{randall}, but these were due to a minor error; in any case, the bounds we report in the present paper improve on the original claims as well.)

\medskip

\noindent {\bf Our second main result}: Here we apply the
combinatorial
insights gained from the proof of Theorem \ref{thm-phasecoexistence} to establish slow mixing of Glauber dynamics
on boxes in $\Z^2$
for values of $\lambda$ that are an order of magnitude lower than the previously best known bounds.
\begin{theorem} \label{thm-mixing}
For all $\lambda > \besttorus$, the mixing time for Glauber dynamics for the hard-core model on $n$ by $n$ boxes in $\Z^2$ with periodic boundary conditions and with activity $\lambda$ is exponential in $n$. For free boundary conditions, we have the same result with $\lambda > \bestbox$.
\end{theorem}

\medskip

\noindent {\bf Our techniques}: The standard approach to showing multiple Gibbs distributions in a statistical physics model is to consider the limiting distributions corresponding to two different boundary conditions on boxes in the lattice centered at the origin,
and find a statistic that separates these two limits. For the hard-core model, it suffices (see \cite{bs}) to compare the {\em even boundary condition} -- all vertices on the boundary of a box at an even distance from the origin are occupied -- and its counterpart the odd boundary condition, and the distinguishing statistic is typically the occupation of the origin. Under odd boundary condition the origin should be unlikely to be occupied, since independent sets with odd boundary and (even) origin occupied must have a {\em contour} -- a two-layer thick unoccupied loop of vertices separating an inner region around the origin that is in ``even phase'' from an outer region near the boundary that is in ``odd phase''. For large enough $\lambda$, such an unoccupied layer is costly, and so such configurations are unlikely. This is essentially the {\em Peierls argument} for phase coexistence, and was the approach taken by Dobrushin \cite{dob}.

As we will see presently, the effectiveness of the Peierls argument is driven by the number of contours of each possible length -- better upper bounds on the number of contours translate directly to better upper bounds on $\lambda_c$. Previous (unpublished) work on phase coexistence in the hard-core model on $\Z^2$ had viewed contours as simple polygons in $\Z^2$, which are closely related to the very well studied family of self-avoiding walks. While this is essentially the best possible point of view when applying the Peierls argument on the Ising model, it is far from optimal for the hard-core model. One of the two major breakthroughs of the present paper is the discovery that hard-core contours, if appropriately defined, can be viewed as simple polygons in the oriented {\em Manhattan lattice} (orient edges of $\Z^2$ that are parallel to the $x$-axis (resp. $y$-axis) positively if their $y$-coordinate (resp. $x$-coordinate) is even, and negatively otherwise), with the additional constraint that contours cannot make two consecutive turns. The number of such polygons can be understood by analyzing a new class of self-avoiding walks, that we refer to as {\em taxi walks}. The number of taxi walks turns out to be significantly smaller than the number of ordinary self-avoiding walks, leading to much better bounds on $\lambda_c$ than could possibly have been obtained previously.

There is a single number $\mu_t>0$, the {\em taxi walk connective constant}, that asymptotically controls the number $\c_n$ of taxi walks of length $n$, in the sense that $\c_n = \mu_t^n f_t(n)$ with $f_t(n)$ sub-exponential. Adapting methods of Alm \cite{Alm} we obtain good estimates on $\mu_t$, giving us good understanding of $\c_n$ for large $n$. The sub-exponential correction makes it difficult to control $\c_n$ for small $n$, however, presenting a major stumbling block to the effectiveness of the Peierls argument. Using the statistic ``occupation of origin'' to distinguish the two boundary conditions, one inevitably has to control $\c_n$ for both small and large $n$. The lack of precise information about the number of short contours leads to discrepancies between theoretical lower limits
and actual bounds, such as that between $C=80$ (theoretical best possible) and $C\approx 300$ for phase coexistence on $\Z^2$ discussed earlier.

The second breakthrough of the present paper is the idea of using an event to distinguish the two boundary conditions that has the property that every independent set in the event has associated with it a {\em long} contour. This allows us to focus exclusively on the asymptotic growth rate of contours/taxi walks, and obviates the need for an analysis of short contours. The immediate result of this breakthrough is that the actual limits of our arguments agree exactly with their theoretical counterparts. The distinguishing event we use extends the idea of {\em fault lines}, which we discuss in more detail below in the context of slow mixing.

The traditional argument for slow mixing is based on the observation that when $\lambda$ is large, the Gibbs distribution favors
dense configurations, and Glauber dynamics will take exponential time to converge to equilibrium.  The slow convergence arises
because the Gibbs distribution is bimodal:  dense configurations lie predominantly  on either the odd or the even sublattice,
while configurations that are roughly half odd and half even have much smaller
probability.  Since Glauber dynamics changes the relative numbers of
even and odd vertices by at most 1 in each step,
the Markov chain has a bottleneck leading to {\em torpid} (slow)  mixing.

Our present work builds on a novel idea from \cite{randall} in which  the notion of {fault lines} was introduced to establish {slow} mixing for the Glauber dynamics on hard-core configurations for moderately large $\lambda$, still improving upon what was best known at that time.
Randall \cite{randall} gave an improvement by realizing that the state space could be partitioned according to certain {\it topological obstructions} in the configurations, rather than the relative numbers of odd or even vertices.  Not only does this approach give better bounds on $\lambda$, but it also greatly simplifies the calculations.  First consider an $n \times n$ lattice region $G$ with free (non-periodic) boundary conditions.  A configuration $I$ is said to have a {\it fault line} if there is a width two path of unoccupied vertices in $I$  from the top of $G$ to the bottom or from the left boundary of $G$ to the right.  Configurations that do not have a fault line must have a {\it  cross} of occupied vertices in either the even or the odd sublattices forming a connected path in $G^2$ from both the top to the bottom and from the left to the right of $G$, where $G^2$ connects vertices at distance 2 in $G$.
Roughly speaking the set of configurations that have a fault line forms a cut set that must be crossed to move from a configuration that has an odd cross to one with an even cross, and it was shown that fault lines are exponentially unlikely when $\lambda$ is large.  Likewise, if $\hat{G}$ is an $n \times n$ region with periodic boundary conditions, it was shown that either there is an odd or an even cross forming non-contractible loops in two different directions or there is a pair of non-contractible fault lines, and a similar argument can be made.

We improve this argument by refining our consideration of fault lines. Previous bounds just characterized fault lines as (rotated) self-avoiding walks in $\Z^2$; here we observe that, suitably modified, they are in fact taxi walks and so the machinery  developed for phase coexistence can be brought to bear in the mixing context, leading to our improvements.

We believe that both of the breakthroughs of the present paper have more general applicability. A natural next step is to study the hard-core model on $\Z^d$ for $d \geq 3$.
Establishing reasonable bounds on the critical activity for $\Z^3$ is an immediate challenge, as is pinning down how the critical value changes with $d$. The best upper bounds are $\widetilde{O}(d^{-1/4})$ for slow mixing \cite{galvin} and $\widetilde{O}(d^{-1/3})$ \cite{peledsamotij} for phase coexistence; the best known lower bounds in both cases are $\Omega(d^{-1})$.

The rest of this manuscript is laid out as follows.
Section~\ref{sec-comb} provides much of the necessary combinatorial background material. This includes the notion of fault lines and crosses, the connection between fault lines and taxi walks, and the characterization of independent sets on finite regions of $\Z^2$. In Section \ref{sec-phasecoexistence} we explain how fault lines and taxi walks can be used to prove Theorem \ref{thm-phasecoexistence} (phase coexistence). In Section~\ref{sec-mixing} we prove Theorem \ref{thm-mixing} (slow mixing). Finally, in Section~\ref{sec-taxi} we give the detailed analysis of the growth rate of taxi walks.

\section{Combinatorial background} \label{sec-comb}
Here we introduce the notions of crosses, fault lines and taxi walks, which are the key ingredients of the proofs of both Theorem \ref{thm-phasecoexistence} and Theorem \ref{thm-mixing}.

\subsection{Crosses and fault lines}\label{termsec}

We begin by defining some useful graph structures.
Let $G=(V, E)$ be a simply connected region in $\Z^2$, say the $n \times n$ square.  We define the graph $\Gd = (\Vd, \Ed)$ as follows.  The vertices $\Vd$ are associated with the midpoints of edges in $E$.  Vertices $u$ and $v$ in $\Vd$ are connected by an edge in $\Ed$ if and only if they are the midpoints of incident edges in $E$ that are perpendicular.  Notice that $\Gd$ is a region in a smaller Cartesian lattice that has been rotated by 45 degrees.

We will also make use of the {\it even and odd subgraphs} of~$G$.  For $b \in \{0, 1\}$, let $\Gb = (\Vb, \Eb)$ be the graph whose vertex set $V_b \subseteq V$ contains all vertices with parity $b$ (i.e., the sum of their coordinates has parity $b$), with $(u, v) \in \Eb$ if $u$ and $v$ are connected in $G^2$.  We refer to $G_0$ and $G_1$ as the even and odd subgraphs.
The graphs $\Gd, \Ge$ and $\Go$ play a central role in defining the features of
independent sets that determine  distinguishing events in our proof of phase coexistence and the partition of the state space for our
proofs of slow mixing.

Given an independent set $I \in \Omega$, we say that a simple path $p$ in $\Gd$ is {\it spanning} if it extends from the top boundary of $\Gd$ to the
bottom, or from the left boundary to the right, and each vertex in $p$ corresponds to an edge in $G$ such that both endpoints are unoccupied in $I$.  It will be convenient to color the vertices in $\Vd$ along a spanning path using the parity of the vertex to the ``left'' (or ``top'') of the path in $V$.
In particular, recall that each vertex $v \in \Vd$ on the path $p$ bisects an edge $e_v \in E$.  Each edge in $E$ has an odd and an even endpoint, and we color $v$ {\it blue} if the odd vertex in $e_v$ is to the left when the path crosses $v$, and {\it red} otherwise. Every time the color of the vertices along the path changes, we have an {\it alternation point}.
It was shown in \cite{randall} that if an independent set has a spanning path, then it must also have one with zero or one alternation points.  We call this path a {\it fault line}, and we let $\Omega_{{\mathcal F}}$ be the set of independent sets in $\Omega$ that contain at least one fault line.

We say that $I \in \Omega$ has an {\em even  bridge} if there is a path from the
left to the right boundary or from the top to the bottom boundary in $\Ge$ consisting of occupied vertices in $I$.
Similarly, we say it has an {\em odd  bridge} if it traverses $\Go$ in either direction.  We say that $I$ has
a {\em  cross} if it has both left-right and a top-bottom  bridges.

Notice that if an independent set has an even top-bottom  bridge it cannot have an odd left-right bridge,
so if it has a  cross, both of its bridges must have the same parity.  We let $\Omega_0 \subseteq \Omega$ be the set of configurations that contain an even  cross and let $\Omega_1 \subseteq \Omega$ be the set of those with an odd   cross.

We can now partition the state space $\Omega$ into three sets,
with one separating the other two; this partition is critical to the proofs of both Theorem \ref{thm-phasecoexistence} and Theorem \ref{thm-mixing}.
The following lemma was proven in \cite{randall}.
\begin{lemma} \label{lem-partition}
Let $G$ be as above. The set $\Omega$ of independent sets on $G$ can be partitioned into sets $\Omega_{{\mathcal F}}$, $\Omega_0$ and $\Omega_1$, consisting of configurations with a fault line, an even cross or an odd cross. If $I \in \Omega_0$ and $I' \in \Omega_1$ then $|I \triangle I'| > 1$.
\end{lemma}

It will be useful to extend these definitions to the torus as well.  Let $n$ be even, and let $\hat{G}$ be the $n \times n$ toroidal region $\{0, \dots, n-1\} \times \{0, \dots n-1\}$, where $v = (v_1, v_2)$ and $u=(u_1, u_2)$ are connected if
$v_1 = u_1 \pm 1 (\mbox{mod }n)$ and $v_2=u_2$ or $v_2 = u_2 \pm 1 (\mbox{mod }n)$ and $v_1=u_1$.
Let $\hat\Omega$ be the set of independent sets on $\hat{G}$ and let $\hat\pi$ be the Gibbs distribution.  As before,
we consider Glauber dynamics that connect configurations with symmetric difference of size  one.

We define $\hatG_{\Diamond}, \hatG_0$ and $\hatG_1$ as above to represent the graph connecting the midpoints of perpendicular edges (including the boundary edges), and the odd and even subgraphs.  As with $\hatG$, all of these have toroidal boundary conditions.

Given $I \in \hat\Omega$, we say that $I$ has a {\em fault} $F = (F_1, F_2)$
if there are a {\it pair} of vertex-disjoint non-contractible cycles $F_1,  F_2$ in $\hatG_{\Diamond}$ whose vertices correspond to edges in $\hatG$ where both endpoints are unoccupied, and such that the vertices on each cycle are all red or all blue (i.e., the endpoints in $\hatG$ to one side of either cycle all have the same parity). We say that $I$ has a {\em  cross} if it has at least two non-contractible cycles of occupied sites in $I$ with different winding numbers.
The next lemma (from \cite{randall}) 
utilizes faults to partition $\hat\Omega$.

\begin{lemma} \label{lem-partition2}
Let $\hatG$ be as above. The set $\hat\Omega$ of independent sets on $\hatG$ can be partitioned into sets $\hof, \hat\Omega_0, \hat\Omega_1$, consisting of configurations with a fault, an even cross or an odd cross.
If $I \in \hat\Omega_0$ and $I' \in \hat\Omega_1$ then $|I \triangle I'| > 1$.
\end{lemma}

\subsection{Taxi walks}\label{taxi}
The strategy for the proofs of phase coexistence and slow mixing will be to use a Peierls argument to define a map from $\of$ to $\Omega$ that takes configurations with fault lines to ones with exponentially larger weight.  The map is not injective, however, so we need to be careful about how large the pre-image of a configuration can be, and for this it is necessary to get a good bound on the number of fault lines.  In \cite{randall} the number of fault lines was bounded by the number of self-avoiding walks in $\Gd$ (or $\hatG_{\Diamond}$ on the torus).  However, this is a gross overcount because, as we
shall see,  this includes all spanning paths with an arbitrary number of alternation points.  As we shall see, we can get much better bounds on the number of fault lines by only counting a subset of
self-avoiding walks with zero or one alternation points.

To begin formalizing this idea, we put an orientation on the edges of $\Gd$.  Each edge $(u, v) \in \Ed$ corresponds to two edges in $E$ that share a vertex $w \in V$.  We orient the edge ``clockwise'' around $w$ if $w$ is even and ``counterclockwise'' around $w$ if $w$ is odd.  For paths with zero alternation points, all of the edges must be oriented in the same direction (with respect to this edge orientation).
If we rotate $\Gd$ so that the edges are axis aligned, then this simply means that the horizontal (resp. vertical) edges alternate direction according to the parity of the $y$- (resp. $x$-) coordinates, like in many well-known metropolises.

We can now define taxi walks.
Let $\Ztaxi$ be an orientation of $\Z^2$ in which an edge parallel to the $x$-axis (resp. $y$-axis) is oriented in the positive $x$-direction if its $y$-coordinate is even (resp. oriented in the positive $y$-direction if its $x$-coordinate is even), and is oriented in the negative direction otherwise (note that this corresponds exactly to the orientation placed on $\Gd$ above).
It is common to refer to $\Ztaxi$ as the {\em Manhattan lattice}: streets are horizontal, with even numbered streets oriented East and odd numbered streets oriented West, and avenues are vertical, with even numbered avenues oriented North and odd numbered avenues oriented South.
\begin{definition}
A {\em taxi walk} is an oriented walk in $\Ztaxi$ that begins at the origin, never revisits a vertex,
and
never takes two left or two right turns in a row.
\end{definition}
We call these taxi walks because the violation of either restriction during a real taxi ride  would cause suspicion among savvy passengers.

\begin{lemma} \label{lem-minimality}
If an independent set $I$ has a fault line $F$ with no alternations, then it also has a fault line $F'$ so that either $F'$ or $F'^R$ (the reversal of $F'$) is a taxi walk.
\end{lemma}

\noindent {\bf Proof}:
It is straightforward to see that if $I$ has a fault line $F$ with no alternation points, then it must have all of its edges oriented the same way (in $\Gd$) and it must be self-avoiding.  Suppose $F$ is a minimal length fault line
in $I$ without any alternations, and suppose that $F$ has two successive turns.  Because of the parity constraints, the vertices immediately before and after these two turns must both connect edges that are in the same direction, and these five edges can be replaced by  a single edge to form a shorter fault line without any alternations.  This is a contradiction to $F$ being minimal, 
completing the proof. 
\qed

\medskip

The same argument  can be used to show that if $F$ is a fault line with an alternation point, then there is a fault line that is the concatenation of two taxi walks (or the reversals of taxi walks). 
Lemma~\ref{lem-minimality}, and the extension just mentioned, are key ingredients in our proofs of both phase coexistence and slow mixing. They allow us to assume, as we do throughout, that all fault lines we work with are essentially taxi walks.
For phase coexistence we will also need to understand the connection between {\em Peierls contours and taxi walks}.

Given an independent set $I$ in $\Z^2$, let $(I^{\mathcal O})^+$ be the set of odd vertices in $I$ together with their neighbors. Let $R$ be any finite component of $(I^{\mathcal O})^+$, and let $W$ be the unique infinite component of $\Z^2 \setminus R$. Let $C$ be the complement of $W$ (the process of going from $R$ to $C$ is essentially one of ``filling in holes'' in $R$). Finally, let $\gamma$ be the set of edges with one end in $W$ and one in $C$, and write $\gamma_{\Diamond}$ for the subgraph of $\Gd$ induced by $\gamma$.
\begin{lemma} \label{lem-contours_are_taxi_walks}
In $\Gd$, $\gamma_{\Diamond}$ is a directed cycle that does not take two consecutive turns. Consequently, if an edge is removed from $\gamma$, the resulting path in $\Ztaxi$ (suitably translated and rotated so that it starts at the origin) is a taxi walk.
\end{lemma}

\noindent {\bf Proof}:
Because $\gamma$ separates $W$ from its complement, $\gamma_{\Diamond}$ must include a cycle surrounding
a vertex of
$W$, and since $\gamma$ is in fact a {\em minimal} edge cutset
($W$ and $C$ are both connected),
$\gamma_{\Diamond}$ must consist of just this cycle.
To see both that $\gamma_{\Diamond}$ is correctly (i.e. uniformly) oriented in $\Gd$, and that it does not take two consecutive turns, note that if either of these conditions were violated then we must have one of the following: a vertex of $W$ (or $C$) all of whose neighbors are in $C$ (or $W$), or a unit square in $\Z^2$ with both even vertices in $C$ and both odd vertices in $W$ (this is an easy case analysis). All of these situations lead to a 4-cycle in $\gamma_{\Diamond}$, a contradiction since $\gamma_{\Diamond}$ is a cycle whose length is evidently greater than $4$ (in fact it must have length at least $12$).
\qed

\medskip



A critical step in our arguments will be bounding the number of taxi walks.  We start by recalling facts about standard self-avoiding walks which have been studied extensively, although many basic questions remain (see, e.g., \cite{ms}).   It is easy to see that on $\Z^2$, the number $c_n$ of walks of length $n$  grows exponentially with $n$ as $2^n \leq c_n \leq 4 \times 3^{n-1}$, since there are always at most 3 ways to extend a self-avoiding walk of length $n-1$ and walks that only take steps to the right or up can always be extended in 2  ways. Hammersley and Welsh \cite{hw} showed that
$c_n =  \mu^n \exp(O(\sqrt{n}))$,
where $\mu$ is known as the {\em connective constant}. It is believed (and there is considerable experimental and heuristical
evidence to suggest) that $\exp(O(\sqrt{n}))$ here can be replaced by $\Theta(n^{11/32})$.

Let $\c_n$ be the number of taxi walks of length $n$.
It is easy to see that $2^{n/2} < \c_n < 4 \times 3^{n-1}$.  For the upper bound here we use $\c_n < c_n$,  and for the lower bound we observe that if we take two steps at a time in one direction we can always go East or North. With very little extra work, we can get a significantly better upper bound:
\begin{lemma}   
Let $\c_n$ be the number of taxi walks of length~$n$.  Then $\c_n = O((1 + \sqrt{5})/2)^n$.
\end{lemma}

\noindent {\bf Proof}:
At each vertex there are exactly two outgoing edges in $\Ztaxi$.  If we arrive at $v$ from $u$, then one of the outgoing
edges continues the walk in the same direction and the other is a turn.
The two allowable directions are determined
by the parity of the coordinates of $v$, so we can encode each walk as a bitstring $s \in \{0, 1\}^{n-1}$.  If $s_1=0$ then the
walk starts by going East (along a street) and if $s_0 = 1$ the walk starts North along an avenue.   For all $i > 1$, if $s_i = 0$ the walk
continues in the same direction as the previous step, while if $s_i = 1$ then the walk turns in the permissible direction.
In this encoding, the condition forbidding consecutive turns forces $s$ to avoid having two $1$'s in a row, and hence $\c_n \leq f_n = O(\phi^n)$,  where
$f_n$ is the $n$th Fibonacci number and $\phi = (1+\sqrt{5})/2 \approx 1.618$ is the golden ratio.
\qed

\medskip

Using more sophisticated tools, described in more detail in Section \ref{sec-taxi}, we can improve our bounds.
\begin{theorem} \label{thm-best_mut}
Let $\c_n$ be the number of taxi walks of length $n$.  Then $\c_n = O(\bestmuupper)^n$. Moreover, there exists a constant $\mu_t >0$ (the {\em taxi walk connective constant}) such that for all $n$, $\c_n = f_t(n) \mu_t^n$, where $f_t(n)$ grows sub-exponentially. We have $\bestmulower < \mu_t < \bestmuupper$ and $\bestlambdalower < \mu_t^4-1 < \besttorus$. If $\mu > \mu_t$ then for all large $n$, $\c_n < \mu^n$.
\end{theorem}

\section{Proof of Theorem 1.1} \label{sec-phasecoexistence}
Here we establish phase coexistence for the hard-core model on $\Z^2$ for all sufficiently large $\lambda$.
The following stronger statement implies Theorem \ref{thm-phasecoexistence} via Theorem~\ref{thm-best_mut}.
\begin{theorem} \label{thm-phasecoexistence_strong}
The hard-core model on $\Z^2$ with activity $\lambda$ admits multiple Gibbs states for all $\lambda > \mu_t^4 -1$, where $\mu_t$ is the connective constant of taxi walks.
\end{theorem}
We will not review the theory of Gibbs states, contenting ourselves with saying informally that an interpretation of the existence of multiple Gibbs states is that the local behavior of a randomly chosen independent set in a box can be made to depend on a boundary condition imposed on the box, even in the limit as the size of the box grows to infinity. See e.g. \cite{Georgii} for a very general treatment, or \cite{bs} for a treatment specific to the hard-core model on the lattice.

Let $U_n$ be the box $[-n,+n]^2$, and $I^{\rm e}$ the independent set consisting of all even vertices of $\Z^2$. Let ${\mathcal J}_n^{\rm e}$ be the set of independent sets that agree with $I^{\rm e}$ off $U_n$, and $\mu_n^{\rm e}$ the distribution supported on ${\mathcal J}_n^{\rm e}$ in which each set is selected with probability proportional to $\lambda^{|I \cap U_n|}$. Define $\mu_n^{\rm o}$ analogously (with ``even'' everywhere replaced by ``odd''). We will exhibit an event ${\mathcal A}$ that depends only on finitely many vertices, with the property that for all large $n$, $\mu_n^{\rm e}({\mathcal A}) \leq 1/3$ and $\mu_n^{\rm o}({\mathcal A}) \geq 2/3$. This is well known (see e.g. \cite{bs}) to be enough to establish the existence of multiple Gibbs states.

The event ${\mathcal A}$ depends on a parameter $m=m(\lambda)$ whose value will be specified later. Specifically, ${\mathcal A}$ consists of all independent sets in $\Z^2$ whose restriction to $U_m$ contains either an odd cross or a fault line. We will show that $\mu_n^{\rm e}({\mathcal A}) \leq 1/3$ for all sufficiently large $n$; reversing the roles of odd and even throughout, the same argument gives that under $\mu_n^{\rm o}$ the probability of $U_m$ having either an {\em even} cross or a fault line is also at most $1/3$, so that (by Lemma \ref{lem-partition}) $\mu_n^{\rm odd}({\mathcal A}) \geq 2/3$.

Write ${\mathcal A}_n^{\rm e}$ for ${\mathcal A} \cap {\mathcal J}_n^{\rm e}$; note that for all large $n$ we have $\mu_n^{\rm e}({\mathcal A}) =  \mu_n^{\rm e}({\mathcal A}_n^{\rm e})$. To show $\mu_n^{\rm e}({\mathcal A}_n^{\rm e}) \leq 1/3$ we will use the fact that $I \in {\mathcal A}_n^{\rm e}$ is in even phase (predominantly even-occupied) outside $U_n$, but because of either the odd cross or the fault line in $U_m$ it is not in even phase close to $U_m$; so there must be a contour marking the furthest extent of the even phase inside $U_n$. We will modify $I$ inside the contour via a weight-increasing map, showing that an odd cross or fault line is unlikely.

\subsection{The contour and its properties} \label{subsec-contours}

Fix $I \in {\mathcal A}_n^{\rm e}$. If $I$ has an odd cross in $U_m$, we proceed as follows (using the notation from the discussion preceding Lemma~\ref{lem-contours_are_taxi_walks}). Let $R$ be the component of $(I^{\mathcal O})^+$ that includes a particular odd cross. Note that because $I$ agrees with $I^{\rm e}$ off $U_n$, $R$ does not reach the boundary of $U_n$, and so  as in the discussion preceding Lemma \ref{lem-contours_are_taxi_walks}, we can associate to $R$ a cutset $\gamma$ separating it from the boundary of $U_n.$

Notice that $\gamma$ is an edge cutset in $U_n$ separating an interior connected region that meets $U_m$ from an exterior connected region that includes the boundary of $U_n$, with all edges from the interior of $\gamma$ to the exterior going from an unoccupied even vertex to an unoccupied odd vertex.
This implies that $|\gamma|$, the number of edges in $\gamma$, is a multiple of~$4$, specifically four times the difference between the number of even and odd vertices in the interior of $\gamma$. Because the interior includes two points of the odd cross that are at distance at least $2m+1$ from each other in $U_m$, we have a lower bound on $\gamma$ that is linear in $m$; in particular, clearly $|\gamma| \geq m$. Note also that by Lemma \ref{lem-contours_are_taxi_walks}, $\gamma_{\Diamond}$ is a closed taxi walk.

We now come to the heart of the Peierls argument. If we modify $I$ by shifting it by one axis-parallel unit (positively or negatively) in the interior of $\gamma$ and leaving it unchanged elsewhere, then the resulting set is still independent, and
we may augment
it
with any vertex in the interior whose neighbor in the direction opposite to the shift is in the exterior. This is a straightforward verification; see \cite[Lemma 6]{bc7} or \cite[Proposition 2.12]{gk} where this is proved in essentially the same setting. Furthermore, from \cite[Lemma 5]{bc7} each of the four possible shift directions free up exactly $|\gamma|/4$ vertices that can be added to the modified independent set.
%
%

We now describe the contour if $I$ has a fault line in $U_m$. If there happens to be an odd occupied vertex in $U_m$ then we construct $\gamma$ as before, starting with some arbitrary component of $(I^{\mathcal O})^+$ that meets $U_m$ in place of the component of an odd cross. If the resulting $\gamma$ has a fault line in its interior, then $\gamma$ and its associated $\gamma_{\Diamond}$ satisfy all the previously established properties immediately.

Otherwise, choose a fault line, which we can assume by Lemma~\ref{lem-minimality} is a taxi walk or the concatenation of two taxi walks.
Whether it has zero or one alternation points, we can find a path $P=u_1u_1\ldots u_k$ in $\Z^2$ with $k$ linear in $m$, with $u_1$ and $u_k$ both odd, with no two consecutive edges parallel,
and with the midpoints of the edges of the path inducing an alternation-free sub-path of the chosen fault line (essentially we are just taking a long piece of the fault line, on an appropriately chosen side of the alternation point, if there is one). This sub-path $F_1$ is a taxi walk. Next, we find a second path in $\Gd$, disjoint from $F_1$, that always bisects completely unoccupied edges, and that taken together with $F_1$ completely encloses $P$. If there are no occupied odd vertices adjacent to even vertices of $P$, such a path is easy to find: we can shift $F_1$ one unit in an appropriate direction, and close off with an additional edge at each end. If there are some odd occupied vertices adjacent to some even vertices of $P$, then this translate of $F_1$ has to be looped around the corresponding components of $(I^{\mathcal O})^+$. Such a looping is possible because $(I^{\mathcal O})^+$ does not reach the boundary of $U_n$, nor does it enclose the fault line (if it did,
  we would be in the case of the previous paragraph).

This second path we have constructed may not be a taxi walk; however, following the proof of Lemma \ref{lem-minimality}, we see that a {\em minimal} path $F_2$ satisfying the conditions of our constructed path is indeed a taxi walk. We take the concatenation of $F_1$ and $F_2$ to be $\gamma_{\Diamond}$ in this case, and take $\gamma$ to be the set of edges that are bisected by vertices of $\gamma_{\Diamond}$. The contours in this case satisfy all the properties of those in the previous case.
The standard strategies outlined in \cite{bc7} and \cite{gk} can easily be used to derive the properties in this case.
The one difference is that now $\gamma_{\Diamond}$ may not be a closed taxi walk; but at worst it is the concatenation of two taxi walks, both of length linear in $m$ (and certainly it can be arranged that each has length at least $m/2$).

\subsection{The Peierls argument} \label{subsec-Peierls}

For $J \in {\mathcal J}_n^{\rm e}$ set $w(J) = \lambda^{|J \cap U_n|}$. Our aim is to show that $w({\mathcal A}_n^{\rm e})/w({\mathcal J}_n^{\rm e}) \leq 1/3$.
For $I \in {\mathcal A}_n^{\rm e}$, let $\varphi(I)$ be the set of independent sets obtained from $I$ by shifting in the interior parallel to $(1,0)$ and adding all subsets of the $|\gamma|/4$ vertices by which the shifted independent set can be augmented. For $J \in \varphi(I)$, let $S$  denote the set of added vertices.
Define a bipartite graph on partite sets ${\mathcal A}_n^{\rm e}$ and ${\mathcal J}_n^{\rm e}$ by joining $I \in {\mathcal A}_n^{\rm e}$ to $J \in {\mathcal J}_n^{\rm e}$ if $J \in \varphi(I)$. Give edge $IJ$ weight $w(I)\lambda^{|S|} = w(J)$ (where $S$ is the set of vertices added to $I$ to obtain $J$).

The sum of the weights of edges out of those $I \in {\mathcal A}_n^{\rm e}$ with $|\gamma(I)|=4\ell$ is $(1+\lambda)^\ell$ times the sum of the weights of those $I$. For each $J \in {\mathcal J}_n^{\rm e}$, the sum of the weights of edges into $J$ from this set of $I$'s is $w(J)$ times the degree of $J$ to the set. If $f(\ell)$ is a uniform upper bound on this degree, then
\begin{equation} \label{Peierlssum}
\frac{w({\mathcal A}_n^{\rm e})}{w({\mathcal J}_n^{\rm e})} \leq \sum_{\ell \geq m/4} \frac{f(\ell)}{(1+\lambda)^\ell}.
\end{equation}
The lower bound on $\ell$ here is crucial. The standard Peierls argument takes ${\mathcal A}$ to be the event that a fixed vertex is occupied, and the analysis of probabilities associated with this event requires dealing with short contours, leading to much weaker bounds than we are able to obtain.

To control $f(\ell)$, observe that for each $J \in {\mathcal J}_n^{\rm e}$ and contour $\gamma$ of length $4\ell$ there is at most one $I$ with $\gamma(I)=\gamma$ such that $J \in \varphi(I)$ ($I$ can be
reconstructed from $J$ and $\gamma$, since the set $S$ of added vertices can easily be identified; cf. \cite[Section 2.5]{gk}). It follows that we may bound $f(\ell)$ by the number of contours of length $4\ell$ with a vertex of $U_m$ in their interiors.

Fix $\mu>\mu_t$. By the properties of contours we have established,
up to translations of contours
this number is at most the maximum of $\mu^{4\ell}$ and
$\sum_{j+k = 4\ell:~j,k\geq m/2} \mu^j\mu^k = 4\ell \mu^{4\ell}$
(for all large $m$, here using Theorem \ref{thm-best_mut}).
The restriction of $\Gd$ to $U_m$ has at most $4(2m+1)^2 \leq 17m^2$ edges, so there are at most this many translates of any particular contour that can have a vertex of $U_m$ in its interior. If follows that we may bound $f(\ell)$ by $68 m^2 \ell \mu^{4\ell}$ and so the sum in (\ref{Peierlssum}) by $\sum_{\ell \geq m} 68 m^2 \ell (\mu^4/(1+\lambda))^\ell$.
For any fixed $\lambda > \mu^4 -1$, there is an $m$ large enough  so that this sum is at most $1/3$; we take any such $m$ to be $m(\lambda)$, completing the proof of phase coexistence.


\section{Proof of Theorem 1.2}\label{sec-mixing}
We now bound the  mixing time of Glauber dynamics by showing that the conductance is exponentially small.

\subsection{Glauber dynamics and mixing time}
Let $G \subset \Z^2$ be an $n\times n$ lattice region
and let $\Omega$ be the set of independent sets on $G$.
Our goal is to sample from $\Omega$ according to the Gibbs distribution,
where each $I\in \Omega$ is assigned probability
$$\pi(I)={\lambda^{|I|}}/{Z},$$
and $Z={\sum_{I'\in\Omega} \lambda^{|I'|}}$ is the normalizing constant.

{\it Glauber dynamics} is a local Markov chain that connects two
independent sets if they have symmetric difference of size one.  The Metropolis probabilities \cite{met} that
force the chain to converge to the Gibbs distribution are given by
$$P(I,I') \ =\  \begin{cases} \frac{1}{2n} \min\left( 1, \lambda^{|I'|-|I|}\right),
& \hbox{if \ } I\oplus I'=1,\cr
1-\sum_{J\sim I} P(I,J), & \hbox{if \ } I=I',\cr
0, & \hbox{otherwise. \ }\cr \end{cases}$$

The conductance, introduced by Jerrum and Sinclair \cite{sj}, is a  good measure of the mixing rate of a chain.
Let
$$\Phi = \min_{S\in \Omega: \pi(S)\leq 1/2} \frac{\sum_{x\in S, y\notin S}
\pi(x)P(x,y)}{\pi(S)},$$
where $\pi(S)=\sum_{x\in S} \pi(x)$ is the weight of the cutset $S$.
The following classical theorem provides the connection between low conductance
and slow mixing.
\begin{theorem} {\rm \cite{sj}} \label{condthm}
For any Markov chain with conductance $\Phi$ we have
$ \frac{\Phi^2}{2} \leq Gap(P) \leq 2 \Phi,$
where $Gap(P)$ is the spectral gap of the transition matrix.

\end{theorem}
\noindent The spectral gap is well-known
to be a measure of the mixing rate of a Markov chain (see, e.g., \cite{sinc}),
so an exponentially small conductance is sufficient to show slow mixing.
Using Theorem \ref{condthm}, our goal is therefore to define a partition that has exponentially small
conductance. The machinery developed in Section \ref{termsec} (in particular Lemmas \ref{lem-partition} and  \ref{lem-partition2}) provide exactly such a partition.
%
%

\subsection{Glauber dynamics on the 2-d torus}\label{slowtorussec}
We are now ready to complete the proof of slow mixing, starting first with the two-dimensional torus.
Let $n$ be even, and let $\hat{G} =\{0, \dots, n-1\} \times \{0, \dots, n-1\}$ be  the $n \times n$ lattice region
with toroidal boundary conditions.   We take $\hat\Omega$ to be the set of
independent sets on $\hatG$ and let $\hat\pi$  be the Gibbs distribution.
Lemma \ref{lem-partition2} shows that $\hat\Omega$ may be partitioned into three sets, $\hat\Omega_{{\mathcal F}}$ (independent sets containing a fault), $\hat\Omega_0$ (independent sets with an even  cross) and $\hat\Omega_1$ (independent sets with an odd cross), and that furthermore
$\hat\Omega_0$ and $\hat\Omega_1$ are not directly connected by moves in the chain $P$.  It remains to show that $\hat\pi(\hof)$ is exponentially smaller than both $\hat\pi(\hat\Omega_0)$ and $\hat\pi(\hat\Omega_1)$.
(Clearly $\hat\pi(\hat\Omega_0)=\hat\pi(\hat\Omega_1)$ by symmetry.)
Notice that on the torus we may assume that fault lines have no alternation points; since they start and end at the same place, the number of alternation points must be even.

Given an independent set $I \in \hof$ with  fault   $F = (F_1, F_2)$,
we partition the vertices of $I$ into two sets, $I_A$ and $I_B,$  
depending on which side of $F_1$ and $F_2$ they lie.
Define the length of a fault to be the total number of edges on this path (or paths) in $\Gd$.  Notice that all fault lines with zero
alternation points have length $N=n+2\ell,$ for some integer $\ell,$ since they all have the same parity.


Let $I'=\sigma(I,F)$ be the configuration formed by shifting $I_A$
one to the right.
Let $F_1'=\sigma(F_1)$ and $F_2'=\sigma(F_2)$ be the images of the fault under this shift.
We define the points that lie in $F_1 \cap F_1'$ and $F_2 \cap F_2'$
to be the points that fall ``in between'' $F$ and $F':=(F_1',F_2')$.

It will be convenient to order the set of possible fault lines so that given a configuration $I \in \hat\Omega_{{\mathcal F}}$ we can  identify its {\it first} fault.
The following results are modified from  \cite{randall} and rely on the new characterization of faults as taxi walks.

\begin{lemma}\label{torlem}
Let $\hat\Omega_F$ be the configurations in $\hof$ with first fault
$F=(F_1, F_2)$.  Write the length of $F$ as $2n+2\ell$.
Then
$$\pi(\hat\Omega_F) \leq  (1 + \lambda)^{-(n+\ell)}.$$
\end{lemma}

\noindent {\bf Proof}:
We define an injection $\phi_F: \hat\Omega_F \times \{0,1\}^{n+\ell}
\hookrightarrow \Omega$ so that $\hat\pi(\phi_F(I,r))=\hat\pi(I)\lambda^{|r|}$.
The injection is formed by cutting the torus $\hatG$ along $F_1$ and
$F_2$ and shifting one of the two connected pieces in any direction by one
unit.  There will be exactly $n+\ell$ unoccupied points near
$F$ that are guaranteed to have only unoccupied neighbors.  We add a subset of
the vertices in this set to $I$ according to bits that are one in the vector $r$.

Given this map, we have
\begin{eqnarray*}
1 & = & \hat\pi(\hat\Omega) \\
  & \geq & \sum_{I \in \hat\Omega_{F}} \ \sum_{r \in \{0,1\}^{n+\ell}} \hat\pi(\phi_{F}(I,r)) \\
  & = & \sum_{I \in \hat\Omega_{F}}  \hat\pi(I) \sum_{r \in \{0,1\}^{n+\ell}} \lambda^{|r|}.
\end{eqnarray*}
\qed

\medskip




\begin{theorem}\label{slowtorusthm}
Let $\hat\Omega$ be the set of independent sets on $\hatG$ weighted by
$\hat\pi(I) = \lambda^{|I|}/Z$, where  $Z=\sum_{I \in \hat\Omega}
\lambda^{|I|}$.
Let $\of$ be the set of independent sets on $\hatG$ with a fault.
Then for any $\lambda > \mu_t^4-1$, there is a constant $c > 0$ such that
$$\hat\pi(\of) \leq e^{-c n}.$$
\end{theorem}

\noindent {\bf Proof}:
Fix $\mu$ satisfying $\lambda > \mu^4 -1 > \mu_t^4 -1$, where $\mu_t$ is the taxi walk connective constant.
Summing over possible locations for the two faults $F_1$ and $F_2$ and
using
Lemma~\ref{torlem},
\begin{eqnarray*}
\hat\pi(\hof) & = & \sum_{F} \hat\pi(\hat\Omega_{F}) \\
&  \leq & \sum_{F}  (1+\lambda)^{-(n+\ell)} \\
&  \leq &  \sum_{i=0}^{(n^2-2n)/2} {n\choose{2}} \mu^{4n+4i} (1+\lambda)^{-(n+i)} \\
& <  & n^2 \sum_{i} \left(\frac{\mu^4}{1+\lambda}\right)^{n+i}.
\end{eqnarray*}
The second inequality here uses Theorem \ref{thm-best_mut}. By our choice of $\mu$ we get (for large $n$)
$\pi(\of) \leq e^{-c n}$ for some constant $c >0$; and we can easily modify this constant to deal with all smaller values of $n$.  \qed

\medskip

From Section~\ref{taxi} we know that $\mu_t^4-1 < \besttorus$.
Combining Theorems~\ref{condthm} and~\ref{slowtorusthm}, we thus get the first part of Theorem \ref{thm-mixing}, as well as the following
stronger result.
\begin{corollary}\label{slowcor}
Fix $\lambda > \mu_t^4-1$. Glauber dynamics for sampling independent sets on the $n\times n$ torus $\hatG$
takes time at least $e^{c n}$ to mix, for some constant $c>0$ (depending on~$\lambda$).
\end{corollary}

\noindent {\bf Proof}:
We will bound the conductance by considering the cut $S=\hat\Omega_0$.
It is clear that $\hat\pi(S) \leq 1/2$ since $\overline{S} = \hat\Omega_{{\mathcal F}} \cup \hat\Omega_1$ and
$\hat\pi(\hat\Omega_0)  = \hat\pi(\hat\Omega_1)$.  Thus,
\begin{eqnarray*}
\Phi &  \leq &  \Phi_S \\
& = & \frac{\sum_{s \in \hat\Omega_0, t \in \hat\Omega_{{\mathcal F}}} \hat\pi(s)P(s,t)}{\hat\pi(\hat\Omega_0)} \\
& = & \frac{\sum_{s \in \hat\Omega_0, t \in \hat\Omega_{{\mathcal F}}} \hat\pi(t)P(t,s)}{\hat\pi(\hat\Omega_0)} \\
& \leq & \frac{\sum_{t \in \hat\Omega_{{\mathcal F}}} \hat\pi(t)}{\hat\pi(\hat\Omega_0)}  \\
& = & \frac{\hat\pi(\hat\Omega_{{\mathcal F}})}{\hat\pi(\hat\Omega_0)}.
\end{eqnarray*}
Given Theorem~\ref{slowtorusthm}, it is trivial to show that $\hat\pi(\hat\Omega_0) > 1/3$, thereby establishing
that the conductance is exponentially small.
It follows from Theorem~\ref{condthm} that Glauber dynamics takes exponential time to converge.
\qed

\subsection{Non-periodic boundary conditions}

For regions with non-periodic boundary conditions we also employ a weight-increasing map from configurations with fault lines by performing a shift and adding vertices. In this setting, however, we not only have to reconstruct the position of the fault line, but we must also encode the part of the configuration lost by the shift due to the finite boundary. In this section we give the proof of the following result, which implies the second part of Theorem \ref{thm-mixing}.
\begin{theorem}\label{gridpf}
Fix $\lambda$ satisfying
\begin{enumerate}
\item $1+\lambda > \mu_t^4$ \ \ and
\item $2(1+\lambda) > \mu_t^2 (1+\sqrt{1+4\lambda})$,
\end{enumerate}
where $\mu_t$ is the taxi walk connective constant.
Glauber dynamics for independent sets on the $n\times n$ grid $G$
takes time at least $e^{c n}$ to mix, for some constant $c=c(\lambda)>0$.
\end{theorem}
Simple algebra reveals that the second condition is satisfied whenever $\lambda^2 + (2-\mu_t^2-\mu_t^4)\lambda +(1-\mu_t^2) > 0$.
Using $\mu_t < \bestmuupper$ we find that both of these conditions are met when $\lambda > \bestbox$, as claimed in Theorem \ref{thm-mixing}.

As before, we partition the state space $\Omega$ into three sets, namely $\of$ (independent sets with a fault line), $\Omega_0$ (independent sets with an even  cross), and $\Omega_1$ (independent sets with an odd cross). From Lemma \ref{lem-partition} we know that these sets do indeed partition $\Omega$, and that furthermore
$\Omega_0$ and $\Omega_1$ are not connected by moves in $P$.  It remains to show that $\pi(\of)$
is exponentially smaller than $\pi(\Omega_0)$ and $\pi(\Omega_1)$.

Let $I \in \of$ be an independent set with vertical fault line $F$.
The fault  $F$ partitions the vertices of $G$ into two sets, $\hbox{Right}(F)$ and $\hbox{Left}(F),$
depending on the side of the fault on which they lie.
Recall that a fault has zero or one alternation points,
and the edges form a path
in $\Gd$.
We will represent length of the path as $N=n+2\ell$, for some $\ell \in \Z$. Strictly speaking this is only correct for fault lines with zero alternation points; but this slight abuse will affect the analysis of what follows by only a constant factor.


Let $I'=\sigma(I,F)$ be the configuration formed by shifting $\hbox{Right}(F)$
one to the right.
We will not be concerned right now
if some vertices ``fall off'' the right side of the region $G$.
Let $F'=\sigma(F)$ be $F$ shifted one to the right.
We define the points that lie in $\hbox{Right}(F) \cap \hbox{Left}(F')$
to be the points that fall ``in between'' $F$ and $F'$.

We use the following result from \cite{randall}.
\begin{lemma}\label{faultlem}
Let $I$ be an independent set with a fault line~$F$.
Let $I'=\sigma(F,I)$  and $F'=\sigma(F)$ be defined as above.
\begin{enumerate}
\item $F$ and $F'$ are both fault lines in $I'$.
\item If we form $I''$ by adding all the points that lie in between
$F$ and $F'$ to $I'$ (except the unique odd point incident to the
alternation point,
if it exists), then $I''$ will be an independent set.
\item If $|F| = n+2\ell$, then there are exactly $n+\ell$ points
that lie in between $F$ and $F'$.
\end{enumerate}
Moreover, $F$ and $F'$ are taxi walks or the concatenation of two taxi walks at an alternation point.
\end{lemma}

Let $I \in \of$ be an independent set with a fault line, which we assume
is vertical.  (If $I$ only has horizontal fault lines, we can rotate
$G$ so that it is vertical; the net effect of ignoring these
independent sets is at most a factor of 2 in the upper bound on $\pi(\of)$, and this will
get incorporated into other constant factors.)
Let $F=F(I)$ be the leftmost fault line.  Let the length
of $F$ be $n+2\ell$, for some integer~$\ell$.

Let $G_{1,n}$ be the $1\times n$ lattice representing the last column of
$G$, and let $J$ be any independent set on $G_{1,n}$.
We further partition $\of$, into
$\cup_{F,J}
\Omega_{F,J}$,
where $I \in \Omega_{F,J}$ if it has leftmost fault line $F$ and is equal
to $J$ when restricted to the last column $G_{1,n}$.

\begin{lemma} Let $F$ be a fault in $G$ with length $n+2\ell$ and
let $\delta$ equal the number of alternation points on $F$
(so $\delta = 0$ or 1).
Let $J$ be an independent set on $G_{1,n}$.  With $\Omega_{F,J}$ defined as
above, we have
$$\pi(\Omega_{F,J}) \leq \lambda^{|J|} (1 + \lambda)^{-(n+\ell-\delta)}.$$
\end{lemma}

\noindent {\bf Proof}:
Let $r\in \{0,1\}^{n+\ell-\delta}$ be any binary vector of length $n + \ell - \delta$ and
let $|r|$ denote the number of bits set to 1, where $|r| \leq n+\ell$.
The main step is to define an injective map $\phi_{F,J}:
\of \times \{0,1\}^{n+\ell} \rightarrow \Omega$ such that, for any $I \in \of$,
$$\pi(\phi_{F,J}(I,r))= \pi(I) \lambda^{-|J|+|r|}.$$
Given this map, we have
\begin{eqnarray*}
1 & = & \pi(\Omega) \\
&  \geq & \sum_{I \in \Omega_{F,J}} \ \sum_{r \in \{0,1\}^{n+\ell-\delta}} \pi(\phi_{F,J}(I,r)) \\
& = & \sum_{I \in \Omega_{F,J}} \ \sum_{r \in \{0,1\}^{n+\ell-\delta}} \pi(I) \lambda^{-|J|+|r|} \\
& = & \sum_{I \in \Omega_{F,J}}  \pi(I) \lambda^{-|J|} \sum_{r \in \{0,1\}^{n+\ell-\delta}} \lambda^{|r|} \\
& = & \sum_{I \in \Omega_{F,J}} \pi(I) \lambda^{-|J|} (1+\lambda)^{n+\ell-\delta} \\
& = & \lambda^{-|J|}(1+\lambda)^{n+\ell-\delta} \ \pi(\Omega_{F,J}).
\end{eqnarray*}

We define the injective map $\phi_{F,J}$ in stages.
For any $I \in \Omega_{F,J}$, we delete the last column (which is equal to $J$).  Next,
recalling that any fault line partitions $G$ into two pieces, we identify
all points in $I$ that fall on the right half and shift these to the right
by one using the map $\sigma(I,F)$.
From Lemma~\ref{faultlem} we know that the number
of points that fall between these two fault lines is $n+\ell$, where $n+2\ell$
is the length of the fault.  The final
step defining the map is to insert new points into the independent set
along this strip between the two faults using the vector $r$, thereby adding
$|r|$ new points.  The new independent set $\phi_{F,J}(I,r)$ has
$|I|-|J|+|r|$ points, and hence has weight $\pi(I)\lambda^{-|J|+|r|}$
\qed

\medskip

\begin{lemma}\label{Jlem}
Let $G_{1,n}$ be a $1 \times n$ strip, and let $\Omega_n$
be the set of independent sets in $G_{1,n}$.  Then
$$\sum_{J\in\Omega_n} \lambda^{|J|} \leq
c \left(\frac{1+\sqrt{1+4\lambda}}{2}\right)^n,$$
for some constant $c$.
\end{lemma}
\noindent {\bf Proof}:
Let $\Omega_i$ be the set of independent sets in $G_{1,i}$
and let $T_i = \sum_{J \in \Omega_i} \lambda^{|J|}.$  Then $T_0=1, T_1=1+\lambda,$
and
$$T_i = T_{i-1} + \lambda T_{i-2}.$$
Solving this Fibonacci-like recurrence yields the lemma.
\qed

\medskip

\begin{theorem}\label{slowthm}
Let $\Omega$ be the set of independent sets on the $n \times n$ lattice
$G$ weighted by
$\pi(I) = \lambda^{|I|}/Z$, where $Z=\sum_{I \in \Omega} \lambda^{|I|}$ is the
normalizing constant.
Let $\of$ be the set of independent sets on $G$ with a fault line.
Then $$\pi(\of) \leq p(n) \ e^{-c' n},$$
for some polynomial $p(n)$ and constant $c'>0,$
whenever $\lambda$ satisfies the hypothesis of Theorem \ref{gridpf}.
\end{theorem}
\noindent {\bf Proof}:
We will make use of the injective map $\phi_{F,J}: \Omega_{F,J} \times
\{0,1\}^N \rightarrow \Omega$, where $N=n+2\ell$ is the length of the fault line.
%
Using Lemma~\ref{Jlem} for the third inequality and Theorem \ref{thm-best_mut} for the fourth,
we have
\begin{eqnarray*}
\pi(\of)  & = & \sum_{F,J} \pi(\Omega_{F,J})\\
&\leq &  \sum_{F,J} \lambda^{|J|} (1+\lambda)^{-(n+\ell-\delta)} \\
& \leq & \lambda \sum_F (1+\lambda)^{-(n+\ell)} \sum_{J\in \Omega_r}  \lambda^{|J|} \\
& \leq & \lambda c \sum_F (1+\lambda)^{-(n+\ell)} \left( \frac{1+\sqrt{1+4\lambda}}{2}\right)^n\\
& \leq & \lambda c \sum_{i=0}^{n^2} n \mu^{2(n+2i)} (1+\lambda)^{-(n+i)}\left( \frac{1+\sqrt{1+4 \lambda}}{2}\right)^{\! n} \\
& = & \lambda c n \sum_{i} \left(\frac{\mu^4}{1+\lambda}\right)^{\! i} \! \!
\left(\frac{\mu^2 (1+\sqrt{1+4\lambda})}{2 (1+\lambda)}\right)^{\! n}\!.
\end{eqnarray*}
This means that we will have $\pi(\of) \leq p(n)  e^{-c' n},$ for some
polynomial $p(n),$  if both
$$
(1+\lambda) > \mu^4,~~~2(1+\lambda) > \mu^2 (1+\sqrt{1+4\lambda})
$$
hold. Taking $\mu$ arbitrarily close to $\mu_t$, the result follows.
\qed

\medskip

\noindent Theorem~\ref{gridpf}  follows exactly as Corollary~\ref{slowcor} in Section~\ref{slowtorussec}.
\vskip.1in

\section{Taxi walks: Bounds and Limits}\label{sec-taxi}
We conclude by justifying the upper bound on the number of taxi walks given in Theorem~\ref{thm-best_mut}, as well as providing a lower bound on $\mu_t$.
It is necessary to first establish the submultiplicativity of $\c_n$ (or, equivalently,  the subadditivity of $\log \c_n$).

\begin{lemma}\label{submult}  Let $\c_n$ be the number of taxi walks of length~$n$ and let $1 \leq i \leq n-1$.  Then $\c_n \leq \c_i \ \c_{n-i}$.
\end{lemma}

\noindent {\bf Proof}: As with traditional self-avoiding walks, the key is to recognize that if we split a taxi walk of length $n$  into two pieces,
the resulting pieces are both self-avoiding.  Let $s = s_1, \dots, s_n$ be a taxi walk of length $n$ and let $1 \leq i \leq n-1$.
Then the initial segment of the walk $s_I = s_1, \dots, s_{i+1}$ is a taxi walk of length $i$.
Let $p=(x,y)$ be the $i$th vertex of the walk $s$.  Let $s_F$ be the final $n-i$ steps of the walk $s$ starting at $p$.  We define
$f(s_F)$ by translating the walk so that $f(p)$ is the origin, reflecting horizontally if $p_x$ is odd and reflecting
vertically if $p_y$ is odd.  Notice that this always produces a valid taxi walk of length $n-i$ and the map $f$ is invertible given~$p$.  Therefore $\c_n \leq \c_i \ \c_{n-i}$.
\qed

\medskip

It follows from Lemma~\ref{submult} that $a_n = \log \c_n$ is subadditive, i.e., $a_{n+m} \leq a_n + a_m$.
By Fekete's Lemma (see, e.g., \cite[Lemma 1.2.2]{steele}) we know that $\lim_{n\rightarrow \infty} a_n /n$ exists and
\begin{equation}\label{fekete}
\lim_{n\rightarrow \infty} \frac{a_n}{n} = \inf \frac{a_n}{n}.
\end{equation}
Thus, we can write the number of taxi walks as $\c_n = \mu_t^n f_t(n)$, where  $\mu_t$ is the connective constant associated  with taxi walks and $f_t(n)$ is subexponential in $n$.
%

Subadditivity gives us a strategy for getting a better bound on $\mu_t$.  From (\ref{fekete}) we see that for all $n$, $\log \c_n /n$ is an upper bound for $\log \mu_t$. We exactly enumerated  taxi walks of length $n$, for $n \leq 60$; see \url{http://nd.edu/~dgalvin1/TD/} for this and other data.
Using $c_{60}=2189670407434$ gives a bound of $\mu_t < 1.6058$. Note that exact counts
for larger $n$ will 
immediately improve our bounds on both $\mu_t$ and $\lambda_c$.

The connective constant for ordinary self-avoiding walks has been well studied, and some of the methods used to obtain bounds there can be adapted to deal with taxi walks. In particular, a method of Alm \cite{Alm} is useful. Fix $n > m >0$. Construct a square matrix $A(m,n)$ whose $ij$ entry counts the number of taxi walks of length $n$ that begin with the $i$th taxi walk of length $m$, and end with the $j$th taxi walk of length $m$, for some fixed ordering of the walks of length $m$. To make sense of this, it is necessary to choose, for each $v \in \Ztaxi$, an orientation preserving map $f_v$ of $\Ztaxi$ that sends the origin to $v$; saying that a walk of length $n$ ends with the $j$th walk of length $m$ means that if the length $m$ terminal
segment of the walk is transformed by $f^{-1}_v$ to start at the origin, where $v$ is the first vertex of the terminal segment, then the result is the $j$th walk of length $m$. Then a theorem of Alm \cite{Alm} may be modified to show that $\mu_t$ is bounded above by $\lambda_1(A(m,n))^{1/(n-m)}$, where $\lambda_1$ indicates the largest positive eigenvalue. (Note that when $m=0$ this recovers the subadditivity bound discussed earlier).

We have calculated $A(20,60)$. This is a square matrix of dimension 20114, and a simple symmetry argument reduces the dimension by a factor of $2$. Using MATLAB, we could estimate the largest eigenvalue of this reduced matrix to obtain $\mu_t < \bestmuupper$ and $\mu_t^4-1 < \besttorus$.

A similar strategy can be used to derive lower bounds on $\mu_t$ in order to determine the theoretical limitations of  our approach of characterizing contours by taxi walks.
We have already given the trivial lower bound $\mu_t \geq \sqrt{2}$. To improve this, we consider {\em bridges} (introduced for ordinary self-avoiding walks by Kesten \cite{Kesten}). A {\em bridge}, for our purposes, is a taxi walk that begins by moving from the origin $(0,0)$ to the point $(1,0)$, never revisits the $y$-axis, and ends by taking a step parallel to the $x$-axis to a point on the walk that has maximum $x$-coordinate over all points in the walk (but note that this maximum does not have to be uniquely achieved at the final point).

Let $b_n$ be the number of bridges of length $n$.  Then bridges are supermultiplicative, i.e., $b_n \geq b_i  b_{n-i}$ (and $\log b_n$ is superadditive).  To see this, note that if $\beta_1$ and $\beta_2$ are bridges, then they both begin and end at vertices
whose $y$-coordinates are even
because they are taking steps to the East.
If the parities of the $x$-coordinates of the first vertices in $\beta_1$ and $\beta_2$  agree, then the concatenation of $\beta_1$ and an appropriate translation of $\beta_2$ is also a bridge; if the parities are different then concatenation of $\beta_1$ with a translation of $\beta_2$ after reflecting horizontally will be a valid bridge.
Notice that the parity of the $x$-coordinate of the two pieces allows us to recover whether a reflection was necessary to keep the walk on the directed Manhattan lattice, so bridges are indeed supermultiplicative.
%
It similarly follows that there are at least $b_n^k$ taxi walks of length $kn$ (just concatenate $k$ length $n$ bridges), so that
$$
\mu_t = \lim_{m \rightarrow \infty} \c_m^{1/m} \geq \lim_{k \rightarrow \infty} (b_n^k)^{1/nk} = b_n^{1/n}.
$$
We have enumerated bridges of length up to $60$, in particular discovering that $b_{60}=80312795498$, leading to $\mu_t > \bestmulower$ and $\mu_t^4-1 > \bestlambdalower$.

A consequence of our lower bound on $\mu_t$ is that our present approach to phase coexistence cannot give anything better than $\lambda_c \leq \bestlambdalower$; this tells us that new ideas will be needed to reach the value of $3.796$ suggested by computations as the true value of $\lambda_c$.


%
%
%

%

\end{document}